# Fractional linear periodic recursions of orders two and three.

## Evgeni Lozitsky

**Abstract.** We study fractional linear recursions of the form $z_n = (a_2 \cdot z_{n-2} + a_1 \cdot z_{n-1} + a_0)/(b_2 \cdot z_{n-2} + b_0)$ and $z_n = (a_3 \cdot z_{n-3} + a_2 \cdot z_{n-2} + a_1 \cdot z_{n-1} + a_0)/(b_3 \cdot z_{n-3} + b_0)$, and find periodic recursions with periods of eight and twelve, which, apparently, were not known before.

## I. Introduction

In the remarkable article [1] study fractional linear recursions of the form:

$$z_n = \frac{a_k \cdot z_{n-k} + a_{k-1} \cdot z_{n-(k-1)} + \ldots + a_i \cdot z_{n-i} + \ldots + a_1 \cdot z_{n-1} + a_0}{z_{n-(k+1)}}.$$

An important feature of such recursions is that the numerator does not include the term $a_{k+1} \cdot z_{n-(k+1)}$. This feature is important because the presence of such terms allows for finding new recursions. Recursions $R'(z_{n-k}, \ldots, z_{n-1})$ and $R(z_{n-k}, \ldots, z_{n-1})$ are considered equivalent if it is true for some a, that $R'(z_{n-k}, \ldots, z_{n-1}) = R(a \cdot z_{n-k}, \ldots, a \cdot z_{n-1})/a$. The result is five recursions, which are actually five one-parameter families:

$$z_n = \frac{1}{z_{n-1}}, \quad z_n = \frac{z_{n-1}}{z_{n-2}}, \quad z_n = \frac{1 + z_{n-1}}{z_{n-2}}, \quad z_n = \frac{1 + z_{n-1} + z_{n-2}}{z_{n-3}}, \quad z_n = \frac{-1 - z_{n-1} + z_{n-2}}{z_{n-3}}. \quad (1)$$

The article [1] is important because it explores the general case, the case of arbitrary n and k.

We will be studying two more specific cases, recursions of order two and three, but not for all possible n:

$$z_n = \frac{a_2 \cdot z_{n-2} + a_1 \cdot z_{n-1} + a_0}{b_2 \cdot z_{n-2} + b_0}, \quad z_n = \frac{a_3 \cdot z_{n-3} + a_2 \cdot z_{n-2} + a_1 \cdot z_{n-1} + a_0}{b_3 \cdot z_{n-3} + b_0}.$$

We will find four recursions from the ones listed above and, in addition, nine more recursions. Two recursions of order two with a period of eight. Four recursions of order two with a period of twelve. And finally, three recursions of order three with a period of twelve. Furthermore, we consider another equivalence relation, which leads to the emergence of two-parameter families.

**Question.** Can the ideas and methods of the article [1] be used in the case of arbitrary n and k for recursions of the form:

$$z_n = \frac{a_{k+1} \cdot z_{n-(k+1)} + a_k \cdot z_{n-k} + \ldots + a_i \cdot z_{n-i} + \ldots + a_1 \cdot z_{n-1} + a_0}{b_{k+1} \cdot z_{n-(k+1)} + b_0} ?$$

In the article [2] George Spahn and Doron Zeilbrger formulate a wish: «It would be very interesting to discover such rational dierence equations with higher periods, that do not trivially follow from the known ones by 'merging'.» It was precisely this wish that motivated me to engage in the study of periodic recursions.

**Acknowledgements.** I am very grateful to Doron Zeilberger for his interest in my research, for his emotional support, for publishing my letters on his web page, and for convincing me to write this article.

## II. Recursions of order two.

Let's consider an arbitrary second-order fractional-linear recursion $z_n = R(z_{n-2}, z_{n-1})$:

$$z_n = \frac{a_2 \cdot z_{n-2} + a_1 \cdot z_{n-1} + a_0}{b_2 \cdot z_{n-2} + b_1 \cdot z_{n-1} + b_0} . \tag{2}$$

Let's also consider the opposite recursion $z_{n-2} = R^{-1}(z_n, z_{n-1})$:

$$z_{n-2} = \frac{b_1 \cdot z_n \cdot z_{n-1} + b_0 \cdot z_n - a_1 \cdot z_{n-1} - a_0}{-b_2 \cdot z_n + a_2} . \tag{3}$$

Next we will consider a particular case where the inverse recursion is also a fractional-linear one. This means that $b_1 = 0$. In this case the original and reverse recursions take the form:

$$z_n = \frac{a_2 \cdot z_{n-2} + a_1 \cdot z_{n-1} + a_0}{b_2 \cdot z_{n-2} + b_0} , \quad z_{n-2} = \frac{b_0 \cdot z_n - a_1 \cdot z_{n-1} - a_0}{-b_2 \cdot z_n + a_2} .$$

**Question:** Is it true that all second-order fractional linear periodic recursions are such?

Also note that $b_2 \neq 0$, otherwise the recursion becomes linear. And finally $a_1 \neq 0$, otherwise the recursion is a trivial 'merger' of two recursions of order one.

Let's define a certain equivalence relation. Let's assume that $G(x) = a \cdot x + b$ – linear transformation, then $G^{-1}(x) = x/a - b/a$ – inverse transformation. Two recursions $R'(z_{n-2}, z_{n-1})$ and $R(z_{n-2}, z_{n-1})$ are equivalent if there exist a and b such that $R'(z_{n-2}, z_{n-1}) = G^{-1}(R(G(z_{n-2}), G(z_{n-1})))$.

After the transformation $a = a_1/b_2$, $b = -b_0/b_2$, the original (2) and reverse (3) recursions will take the simplest form:

$$z_n = R(z_{n-2}, z_{n-1}) = \frac{a_2 \cdot z_{n-2} + z_{n-1} + a_0}{z_{n-2}} , \quad z_{n-2} = R^{-1}(z_n, z_{n-1}) = \frac{z_{n-1} + a_0}{z_n - a_2} .$$

Let's now consider two recursions $z_1, z_2, z_3 = R(z_1, z_2), \ldots$ and $u_1, u_2, u_3 = R^{-1}(u_1, u_2), \ldots$ If the original recursion is periodic with a period of k, and if $u_1 = z_2$, $u_2 = z_1$, then $u_j = z_i$, if $i+j = k+3$. The equalities $u_j = z_i$ actually serve as conditions for $a_0$ and $a_2$. The elements $z_i$ and $u_j$ are rational functions of the variables $z_1$ and $z_2$. The equality $z_i = u_j$ must hold identically in terms of $z_1$ and $z_2$, leading to a system of polynomial equations on $a_2$ and $a_0$.

| i+j=k+3 | | | | | | |
|---------|---------|---------|-----------|---------|-----------|-----------|
| $z_1$ | $z_2$ | $z_3$ | .......... $z_i$ .......... | $z_k$ | $z_{k+1}$ | $z_{k+2}$ |
| $\|\|$ | $\|\|$ | $\|\|$ | $\|\|$ | $\|\|$ | $\|\|$ | $\|\|$ |
| $u_{k+2}$ | $u_{k+1}$ | $u_k$ | .......... $u_j$ .......... | $u_3$ | $u_2$ | $u_1$ |

## II_I. Period five.

The main equality is: $z_4 = u_4$. Here is the Wolfram Mathematica code:

```
z3=(a2*z1+z2+a0)/z1
z4=Factor[(a2*z2+z3+a0)/z2]
u1=z2
u2=z1
u3=(u2+a0)/(u1-a2)
u4=Factor[(u3+a0)/(u2-a2)]
CoefficientList[Numerator[Factor[u4-z4]],{z1,z2}]
```

After factoring these coefficients, we obtain a system of equations:

$$a_2=0 \qquad (a_0-a_2)\cdot a_2=0 \qquad -1+a_0+a_2^2=0 \qquad a_2\cdot(-1-a_2+a_2^2)=0$$
$$a_0\cdot a_2^2=0 \qquad a_2\cdot(a_0+a_2)=0 \qquad 1-a_0-a_2+a_2^2=0 \qquad a_2\cdot(-a_0+a_0\cdot a_2+a_2^2)=0$$

This system has the unique solution: $\{a_0=1, a_2=0\}$. This is known as the Lynnes cycle, one of the recursions in the list (1).

## II_II. Period six.

The main equality is: $z_5=u_4$. The Wolfram Mathematica code is analogous.
After factoring the coefficients, we obtain a system of equations:

$$a_2=0 \qquad\qquad a_2\cdot(a_0+2\cdot a_2+a_2^2)=0$$
$$a_0\cdot a_2^2=0 \qquad\qquad a_2\cdot(-a_0+a_0\cdot a_2+a_2^2)=0$$
$$a_2\cdot(a_0+a_2)=0 \qquad\qquad -a_0-a_0^2-a_2+a_0\cdot a_2-a_0\cdot a_2^2+a_2^3=0$$
$$a_0+a_2+a_2^2=0 \qquad\qquad -a_0^2-a_0\cdot a_2+a_0^2\cdot a_2+a_2^2+a_0\cdot a_2^3=0$$
$$a_0-a_0\cdot a_2-a_2^2-a_2^3=0 \qquad\qquad a_2\cdot(-1-2\cdot a_0-a_2+a_0\cdot a_2+a_2^2+a_2^3)=0$$

This system has, obviously, the unique solution: $\{a_0=0, a_2=0\}$. And once again, we obtain the recursion from the list (1).

## II_III. Period eight.

The main equality is: $z_6=u_5$. The Wolfram Mathematica code is similar. After factoring the coefficients, we obtain a system of equations.

$$a_0-a_0\cdot a_2-a_2^3=0 \qquad -a_0-a_2+a_0\cdot a_2+2\cdot a_0\cdot a_2^2+a_2^3+a_2^4=0$$
$$(a_0+a_2)\cdot(1+a_2^2)=0 \qquad -1-2\cdot a_0-a_0\cdot a_2+a_2^2-a_0\cdot a_2^2-2\cdot a_0\cdot a_2^3+a_2^4-a_2^5=0$$
$$a_2\cdot(2\cdot a_0+a_2+a_2^2)=0 \qquad a_0+2\cdot a_0^2-a_0\cdot a_2-2\cdot a_0^2\cdot a_2+2\cdot a_0\cdot a_2^2-2\cdot a_0\cdot a_2^3+a_2^4=0$$
$$-a_0-a_2-a_0\cdot a_2-a_2^2-a_2^3=0 \qquad a_2\cdot(-1-2\cdot a_0-2\cdot a_0^2-a_2+a_0\cdot a_2-a_2^2+a_0\cdot a_2^2+2\cdot a_2^3+a_2^4)=0$$
$$a_0^2\cdot(a_0-a_0\cdot a_2+2\cdot a_0\cdot a_2^2+a_2^4)=0 \qquad a_0\cdot(-2\cdot a_0-a_0^2+2\cdot a_0\cdot a_2+a_0^2\cdot a_2-4\cdot a_0\cdot a_2^2+a_0\cdot a_2^3-2\cdot a_2^4)=0$$
$$a_0+a_0\cdot a_2+a_2^2-2\cdot a_0\cdot a_2^2-a_2^3-a_2^4=0 \qquad a_2\cdot(3+2\cdot a_0^2-a_2-2\cdot a_0^2\cdot a_2+a_2^2+a_0\cdot a_2^2-a_2^3+a_0\cdot a_2^3+a_2^5)=0$$

$$1+3\cdot a_0^2+a_0\cdot a_2+a_0^2\cdot a_2-2\cdot a_0^2\cdot a_2^2-a_2^3+a_0\cdot a_2^3-a_2^4-a_0\cdot a_2^4+a_2^5=0$$
$$a_2\cdot(-3\cdot a_0-a_0^2+2\cdot a_0^3-a_0\cdot a_2+4\cdot a_0^2\cdot a_2-2\cdot a_0\cdot a_2^2-a_0^2\cdot a_2^2-a_2^3-a_0\cdot a_2^4-a_2^5)=0$$
$$-a_0-a_0^2+a_0\cdot a_2-a_0^2\cdot a_2+2\cdot a_2^2-3\cdot a_0\cdot a_2^2-2\cdot a_0^2\cdot a_2^2-a_2^3-a_0\cdot a_2^3+a_2^4-a_0\cdot a_2^4-a_2^5=0$$
$$a_0\cdot(a_0+a_0^2+3\cdot a_0\cdot a_2-3\cdot a_0^2\cdot a_2-a_0\cdot a_2^2+2\cdot a_0^2\cdot a_2^2+a_2^3+3\cdot a_0\cdot a_2^3+a_0\cdot a_2^4+2\cdot a_2^5)=0$$
$$2\cdot a_0+a_0^2+2\cdot a_0^3+3\cdot a_0\cdot a_2-2\cdot a_0^2\cdot a_2-a_0\cdot a_2^2+a_0^2\cdot a_2^2+a_2^3+2\cdot a_0\cdot a_2^3+2\cdot a_0^2\cdot a_2^3+2\cdot a_0\cdot a_2^5+a_0\cdot a_2^5=0$$

The system is quite cumbersome, but it can be solved in several steps, starting with the simplest equation: $(a_0+a_2)\cdot(1+a_2{}^2)=0$. As a result, we obtain two solutions: $\{a_0=(1-i)/2,\ a_2=i\}$, $\{a_0=(1+i)/2,\ a_2=-i\}$.

$$z_n = \frac{2\cdot i\cdot z_{n-2}+z_{n-1}+(1-i)}{2\cdot z_{n-2}},\quad z_n = \frac{-2\cdot i\cdot z_{n-2}+z_{n-1}+(1+i)}{2\cdot z_{n-2}}.$$

## II_IV. Period twelve.

Main equation: $z_8=u_7$. The Wolfram Mathematica code is similar. The system of equations turns out to be quite cumbersome. However, it is still possible to solve this system using the Solve[,{}] command from the Wolfram Mathematica Language. As a result, four solutions are obtained:

$$\left\{a_0=\left(1+\frac{i}{2}\right)-\frac{\sqrt{3}}{2},a_2=-i\right\},\quad \left\{a_0=\left(1+\frac{i}{2}\right)+\frac{\sqrt{3}}{2},a_2=-i\right\},$$

$$\left\{a_0=\left(1-\frac{i}{2}\right)-\frac{\sqrt{3}}{2},a_2=i\right\},\quad \left\{a_0=\left(1-\frac{i}{2}\right)+\frac{\sqrt{3}}{2},a_2=i\right\}.$$

Note that in each of the four cases, $a_0=1+\{$the sixth root of $-1\}$, $a_2=\{$the square root of $-1\}$.

$$z_n = \frac{-2\cdot i\cdot z_{n-2}+2\cdot z_{n-1}+(2-\sqrt{3}+i)}{2\cdot z_{n-2}},\quad z_n = \frac{-2\cdot i\cdot z_{n-2}+2\cdot z_{n-1}+(2+\sqrt{3}+i)}{2\cdot z_{n-2}},$$

$$z_n = \frac{2\cdot i\cdot z_{n-2}+2\cdot z_{n-1}+(2-\sqrt{3}-i)}{2\cdot z_{n-2}},\quad z_n = \frac{2\cdot i\cdot z_{n-2}+2\cdot z_{n-1}+(2+\sqrt{3}-i)}{2\cdot z_{n-2}}.$$

Recursions with periods of four, seven, nine, ten, and eleven were not found.

## III. Recursions of order three.

Let's consider an arbitrary third-order fractional-linear recursion $z_n=R(z_{n-3},z_{n-2},z_{n-1})$:

$$z_n = \frac{a_3\cdot z_{n-3}+a_2\cdot z_{n-2}+a_1\cdot z_{n-1}+a_0}{b_3\cdot z_{n-3}+b_2\cdot z_{n-2}+b_1\cdot z_{n-1}+b_0}. \tag{4}$$

Let's also consider the opposite recursion $z_{n-3}=R^{-1}(z_n,z_{n-1},z_{n-2})$:

$$z_{n-3} = \frac{b_2\cdot z_n\cdot z_{n-2}+b_1\cdot z_n\cdot z_{n-1}-a_2\cdot z_{n-2}-a_1\cdot z_{n-1}+b_1\cdot z_n-a_1}{-b_3\cdot z_n+a_3}. \tag{5}$$

Next we will consider a particular case where the inverse recursion is also a fractional-linear one. This means that $b_1=0$, and $b_2=0$. In this case the original and reverse recursions take the form:

$$z_n = \frac{a_3\cdot z_{n-3}+a_2\cdot z_{n-2}+a_1\cdot z_{n-1}+a_0}{b_3\cdot z_{n-3}+b_0},\quad z_{n-3} = \frac{-a_2\cdot z_{n-2}-a_1\cdot z_{n-1}+b_0\cdot z_n-a_1}{-b_3\cdot z_n+a_3}.$$

**Question:** Is it true that all third-order fractional-linear periodic recursions are such?

Also note that $b_3 \neq 0$, otherwise the recursion becomes linear. And finally $a_1 \neq 0$ or $a_2 \neq 0$, otherwise the recursion is a trivial 'merger' of two recursions of order one.

Two recursions $R'(z_{n-3}, z_{n-2}, z_{n-1})$ and $R(z_{n-3}, z_{n-2}, z_{n-1})$ are equivalent if there exist $a$ and $b$ such that $R'(z_{n-3}, z_{n-2}, z_{n-1}) = G^{-1}(R(G(z_{n-3}), G(z_{n-2}), G(z_{n-1})))$.

Let's consider two types of recursions. Type 1: $a_2 \neq 0$. Type 2: $a_2 = 0$, $a_1 \neq 0$. After transformations, for the first type, $a = a_2/b_3$, $b = -b_0/b_3$, and for the second type, $a = a_1/b_3$, $b = -b_0/b_3$, the original (4) and inverse (5) recursions take a simpler form:

Type 1: $z_n = R(z_{n-3}, z_{n-2}, z_{n-1})$, $z_{n-3} = R^{-1}(z_n, z_{n-1}, z_{n-2})$, where:

$$R(z_{n-3}, z_{n-2}, z_{n-1}) = \frac{a_3 \cdot z_{n-3} + z_{n-2} + a_1 \cdot z_{n-1} + a_0}{z_{n-3}}, \quad R^{-1}(z_n, z_{n-1}, z_{n-2}) = \frac{z_{n-2} + a_1 \cdot z_{n-1} + a_0}{z_n - a_3}.$$

Type 2: $z_n = R(z_{n-3}, z_{n-2}, z_{n-1})$, $z_{n-3} = R^{-1}(z_n, z_{n-1}, z_{n-2})$, where:

$$R(z_{n-3}, z_{n-2}, z_{n-1}) = \frac{a_3 \cdot z_{n-3} + z_{n-1} + a_0}{z_{n-3}}, \quad R^{-1}(z_n, z_{n-1}, z_{n-2}) = \frac{z_{n-1} + a_0}{z_n - a_3}.$$

Let's now consider two recursions: $z_1$, $z_2$, $z_3$, $z_4 = R(z_1, z_2, z_3)$, ... and $u_1$, $u_2$, $u_3$, $u_4 = R^{-1}(u_1, u_2, u_3)$, ... If the original recursion is periodic with a period of $k$, and if $u_1 = z_3$, $u_2 = z_2$, $u_3 = z_1$ then $u_j = z_i$, if $i + j = k + 4$. The equalities $u_j = z_i$ actually serve as conditions for $a_0$, $a_1$, $a_3$ for recursions of the first type, and for $a_0$, $a_3$ for recursions of the second type. The elements $z_i$ and $u_j$ are rational functions of the variables $z_1$, $z_2$, $z_3$. The equality $z_i = u_j$ must hold identically in terms of $z_1$, $z_2$, $z_3$, leading to a system of polynomial equations on $a_0$, $a_1$, $a_3$ and $a_0$, $a_3$ respectively.

| **i+j=k+4** | | | | | | | | |
|---|---|---|---|---|---|---|---|---|
| $z_1$ | $z_2$ | $z_3$ | $z_4$ | .......... $z_i$ .......... | $z_k$ | $z_{k+1}$ | $z_{k+2}$ | $z_{k+3}$ |
| ‖ | ‖ | ‖ | ‖ | ‖ | ‖ | ‖ | ‖ | ‖ |
| $u_{k+3}$ | $u_{k-2}$ | $u_{k+1}$ | $u_k$ | .......... $u_j$ .......... | $u_4$ | $u_3$ | $u_2$ | $u_1$ |

## III_I. Period eight.

The main equality is: $u_6 = z_6$. Here is the Wolfram Mathematica code for Type 1:

```
z4=Factor[(a3*z1+z2+a1*z3+a0)/z1]
z5=Factor[(a3*z2+z3+a1*z4+a0)/z2]
z6=Factor[(a3*z3+z4+a1*z5+a0)/z3]
u1=z3
u2=z2
u3=z1
u4=Factor[(u3+a1*u2+a0)/(u1-a3)]
u5=Factor[(u4+a1*u3+a0)/(u2-a3)]
u6=Factor[(u5+a1*u4+a0)/(u3-a3)]
CoefficientList[Numerator[Factor[u6-z6]], {z1, z2, z3}]
```

After factoring these coefficients, we obtain a system of equations:

$a_3=0$ $\qquad$ $a_0-a_1+a_3+a_1 \cdot a_3-a_3{}^2=0$

$-a_0+a_1-a_3{}^2=0$ $\qquad$ $a_1{}^2 \cdot a_3{}^2 \cdot (a_0-a_1 \cdot a_3)=0$

$a_0 \cdot a_1{}^2 \cdot a_3{}^3=0$ $\qquad$ $a_3 \cdot (a_0+a_1{}^2-a_3-a_1 \cdot a_3)=0$

$(-1+a_1) \cdot (1+a_1)=0$ $\qquad$ $-a_0+a_1 3+a_0 \cdot a_3-2 \cdot a_1 \cdot a_3+a_3 3=0$

$a_1 \cdot (a_1{}^2-a_3) \cdot a_3=0$ $\qquad$ $a_3{}^2 \cdot (a_0 \cdot a_1{}^2-a_0 \cdot a_3-a_1{}^2 \cdot a_3)=0$

$a_3 \cdot (a_0+a_3+a_1 \cdot a_3)=0$ $\qquad$ $a_1 \cdot a_3{}^2 \cdot (a_0 \cdot a_1-a_0 \cdot a_3-a_1 \cdot a_3{}^2)=0$

$(a_0+a_1{}^2-a_3) \cdot a_3{}^2=0$ $\qquad$ $a_3 \cdot (-a_0-a_1{}^2+a_3+a_0 \cdot a_3+a_3{}^2+a_1{}^2 \cdot a_3{}^2)=0$

$a_1 \cdot a_3{}^2 \cdot (a_0+a_1 \cdot a_3)=0$ $\qquad$ $a_3 \cdot (a_0 \cdot a_1-a_0 \cdot a_3+a_1{}^2 \cdot a_3-a_3{}^2-a_1 \cdot a_3{}^2)=0$

$a_1 \cdot a_3 \cdot (a_0-a_3+a_1 \cdot a_3)=0$ $\qquad$ $a_1 \cdot a_3 \cdot (a_0 \cdot a_1-a_0 \cdot a_3-a_1{}^2 \cdot a_3+a_3{}^2-a_1 \cdot a_3{}^2)=0$

$a_3 \cdot (a_0 \cdot a_1{}^2-a_0 \cdot a_3-a_1{}^2 \cdot a_3-a_1{}^3 \cdot a_3+a_1 \cdot a_3{}^2)=0$

$-1+a_0 \cdot a_1-a_0 \cdot a_3+a_1 \cdot a_3+a_1{}^2 \cdot a_3-a_3{}^2-a_1 \cdot a_3{}^2+a_3{}^3=0$

$a_0-a_1-a_0 \cdot a_1+a_1{}^2-a_3-a_1 \cdot a_3+a_1{}^2 \cdot a_3-a_3{}^2-a_1 \cdot a_3{}^2+a_3{}^3=0$

$-a_0+a_0 \cdot a_1{}^2-a_1{}^2 \cdot a_3-a_1{}^3 \cdot a_3+2 \cdot a_1 \cdot a_3{}^2-a_1{}^2 \cdot a_3{}^2+a_3{}^3+a_1 \cdot a_3{}^3-a_3{}^4=0$

$a_3 \cdot (a_0 \cdot a_1{}^2-a_0 \cdot a_3-a_0 \cdot a_1 \cdot a_3-a_1{}^2 \cdot a_3+a_0 \cdot a_3{}^2-a_1{}^2 \cdot a_3{}^2+a_3{}^3+a_1 \cdot a_3{}^3)=0$

This system has, the two solutions: $\{a_0=-1, a_1=-1, a_3=0\}$, $\{a_0=1, a_1=1, a_3=0\}$. And once again, we obtain the recursions from the list (1).

Wolfram Mathematica code for Type 2:

```
z4=Factor[(a3*z1+z3+a0)/z1]
z5=Factor[(a3*z2+z4+a0)/z2]
z6=Factor[(a3*z3+z5+a0)/z3]
u1=z3
u2=z2
u3=z1
u4=Factor[(u2+a0)/(u1-a3)]
u5=Factor[(u3+a0)/(u2-a3)]
u6=Factor[(u4+a0)/(u3-a3)]
CoefficientList[Numerator[Factor[u6-z6]], {z1, z2, z3}]
```

The system of equations has no solutions because one of the coefficients is equal to 1.

## III_II. Period twelve.

The main equality is: $u_8=z_8$. The Wolfram Mathematica code is similar. For recursions of type 1 the system of equations turns out to be very cumbersome, and it has three solutions:

$$\left\{a_0=-\frac{1}{2}, a_1=-1, a_3=-1\right\},$$

$$\left\{a_0=\frac{1}{2}\left(-1-\sqrt{3} \cdot \mathrm{i}\right), a_1=\frac{1}{2}\left(1-\sqrt{3} \cdot \mathrm{i}\right), a_3=\frac{1}{2}\left(1+\sqrt{3} \cdot \mathrm{i}\right)\right\}, \qquad (6)$$

$$\left\{ a_0 = \frac{1}{2}\left(-1+\sqrt{3}\cdot \mathrm{i}\right), a_1 = \frac{1}{2}\left(1+\sqrt{3}\cdot \mathrm{i}\right), a_3 = \frac{1}{2}\left(1-\sqrt{3}\cdot \mathrm{i}\right) \right\}. \qquad (7)$$

The corresponding recursions are:

$$z_n = \frac{-1-2\cdot z_{n-3}+2\cdot z_{n-2}-2\cdot z_{n-1}}{2\cdot z_{n-3}},$$

$$z_n = \frac{\left(1+\sqrt{3}\cdot \mathrm{i}\right)\cdot z_{n-3}+2\cdot z_{n-2}+\left(1-\sqrt{3}\cdot \mathrm{i}\right)\cdot z_{n-1}+(-1-\sqrt{3}\cdot \mathrm{i})}{2\cdot z_{n-3}},$$

$$z_n = \frac{\left(1-\sqrt{3}\cdot \mathrm{i}\right)\cdot z_{n-3}+2\cdot z_{n-2}+\left(1+\sqrt{3}\cdot \mathrm{i}\right)\cdot z_{n-1}+(-1+\sqrt{3}\cdot \mathrm{i})}{2\cdot z_{n-3}}.$$

The first of these three recursions is apparently the most interesting result of this article, as its coefficients are integers.

Note that all the coefficients in (6) and (7) are roots of unity of the sixth degree.

For recursions of type two the system of equations has no solutions, because one of the coefficients is equal to 1.

Recursions with periods of five, six, seven, nine, ten, and eleven were not found.

# References

[1] M. Csöornyei and M. Laczkovich, Some periodic and non-periodic recursions, Monatsh. Math. 132 (2001), 215–236.

[2] George Spahn and Doron Zeilberger, Experimenting with Discrete Dynamical Systems, https://arxiv.org/pdf/2306.11929.pdf, https://sites.math.rutgers.edu/~zeilberg/mamarim/mamarimPDF/dds.pdf